\documentclass[11pt]{article}

\usepackage{amsfonts,amssymb,latexsym,indentfirst}

\newtheorem{prop}{Proposition}
\newtheorem{rema}{Remark}

\newtheorem{lemm}{Lemma}
\newtheorem{theo}{Theorem}
\newtheorem{coro}{Corollary}

\newcommand{\N}[1][]{\ensuremath{{\mathbb{N}^{#1}} }}
\newcommand{\Z}[1][]{\ensuremath{{\mathbb{Z}^{#1}} }}
\newcommand{\C}[1][]{\ensuremath{{\mathbb{C}^{#1}} }}
\newcommand{\R}[1][]{\ensuremath{{\mathbb{R}^{#1}} }}

\def\Re{ \mathrm{Re}\, }

\author{Didier Pilod}
\title{\textbf{On the Cauchy problem for higher-order nonlinear dispersive equations}}
\date{}
\begin{document}
\maketitle

\begin{abstract}
We study the higher-order nonlinear dispersive equation
\begin{displaymath}
\partial_tu+\partial_x^{2j+1}u=\sum_{0\le j_1+j_2 \le 2j}
a_{j_1,j_2}\partial_x^{j_1}u \partial_x^{j_2}u, \quad x, \ t \in \R.
\end{displaymath}
where $u$ is a real- (or complex-) valued function. We show that the
associated initial value problem is well posed in weighted Besov and
Sobolev spaces for small initial data. We also prove ill-posedness
results when $a_{0,k} \neq 0$ for some $k>j$, in the sense that this
equation cannot have its flow map $C^2$ at the origin in $H^s(\R)$,
for any $s \in \R$. The same technique leads to similar
ill-posedness results for other higher-order nonlinear dispersive
equation as higher-order Benjamin-Ono and intermediate long wave
equations.
\end{abstract}

\section{Introduction}
In this paper we consider the initial value problem (IVP)
\begin{equation} \label{sKdV}
\left\{\begin{array}[pos]{ll}
          \partial_tu+\partial_x^{2j+1}u=\sum_{0\le j_1 + j_2 \le 2j}
          a_{j_1,j_2}\partial_x^{j_1}u \partial_x^{j_2}u, \quad x, \ t \in \R, \\
          u(0)=\phi, \\
\end{array} \right.
\end{equation}
where $u$ is a real- (or complex-) valued function and $a_{l_1,l_2}$
are constants in $\R$ or $\C$. It is a particular case of the class
of IVPs
\begin{equation} \label{KdV_h}
\left\{  \begin{array}[pos]{ll}
          \partial_tu+\partial_x^{2j+1}u+P(u,\partial_xu,\ldots,\partial_x^{2j}u),
          \quad x,\ t \in \R, \ j \in \N \\          u(0)=u_0, \\
     \end{array} \right.
\end{equation}
where
\begin{displaymath}
P: \R^{2j+1}\rightarrow \R \quad (\mbox{or} \ P:
\C^{2j+1}\rightarrow \C)
\end{displaymath}
is a polynomial having no constant or linear terms.

The class of IVPs (\ref{KdV_h}) contains the KdV hierarchy as well
as higher-order models in water waves problems (see \cite{KPV4} for
the references). When $j=1$ and the nonlinearity has the form
$u\partial_xu$, the equation (\ref{sKdV}) is the KdV equation, when
$j=1$ and the nonlinearity has the form
$\alpha(\partial_xu)^2+\gamma u\partial^2_xu$, it becomes the limit
(when the dissipation tends to zero) of the KdV-Kuramoto-Velarde
(KdV-KV) equation (see \cite{Arg} and \cite{Pil}).

Kenig, Ponce and Vega have proved that the class of IVPs
(\ref{KdV_h}) is well-posed in some weighted Sobolev spaces for
small initial data \cite{KPV3}, and for arbitrary initial data
\cite{KPV4}. In \cite{Arg} Argento found the best exponents of the
weighted Sobolev spaces where well-posedness for the non dissipative
KdV-KV equation is satisfied. More precisely, she showed that this
IVP is well-posed for small initial data in $H^k(\R) \cap
H^3(\R;x^2dx)$ for $k \in \N$, $k \ge 5$.

The method used, in the case of small initial data, is an
application of a fixed point theorem to the associated integral
equation, taking advantage of the smoothing effects associated to
the unitary group of the linear equation. In particular, a maximal
(in time) function estimate is needed in $L^1_x$. Actually, as
observed in \cite{KPV2}, the $L^1_x$-maximal function estimate fails
without weight. In the case of arbitrary initial data, Kenig, Ponce
and Vega performed a gauge transformation on the equation
(\ref{KdV_h}) to get a dispersive system whose nonlinear terms are
independent of the higher-order derivative. This allows to apply the
techniques already used in the case of small initial data.

In the following, we improved these results for the IVP (\ref{sKdV})
in the case of small initial data, using weighted Besov spaces. The
use of Besov spaces is inspired by the works of Molinet and Ribaud
on the Korteweg-de Vries equation  \cite{MR1} and on the
Benjamin-Ono equation \cite{MR2}, and of Planchon on the nonlinear
Schr\"odinger equation \cite{Pla}. It allows to refine the
$L^1_x$-maximal function estimate, using the $L^4 _x$-maximal
function estimate derived by Kenig and Ruiz \cite{KR} (see also
\cite{KPV1}), and to obtain well-posedness results in fractional
weighted Besov spaces.

Nevertheless, the natural spaces to show well-posedness for the
equation (\ref{sKdV}) are the Sobolev spaces $H^s(\R)$. We prove
here that if there exists $k>j$ such that $a_{0,k} \neq 0$, we
cannot solve this problem in any space continuously embedded in
$C([-T,T];H^s(\R))$, for any $s \in \R$, using a fixed point theorem
on the integral equation. As a consequence of this result, we deduce
that in this case, the flow-map data solution of (\ref{sKdV}) cannot
be $C^2$ at the origin from $H^s(\R)$ to $H^s(\R)$, for any $s \in
\R$.

The same kind of argument leads to similar results for other
higher-order nonlinear dispersive equations. We consider first a
higher order Benjamin-Ono equation.
\begin{equation}\label{ho.BO}
\left\{\begin{array}[pos]{ll}
           \partial_tu-bH\partial^2_xu+
           a\epsilon \partial_x^3u=cu\partial_xu-d\epsilon
           \partial_x(uH\partial_xu+H(u\partial_xu)) \\
           u(0)=\phi, \\
       \end{array} \right.
\end{equation}
where $H$ is the Hilbert transform, $u$ is a real-valued function,
and $a \in \R$, $b, \ c$ and $d$ are positive constants. This
equation was derived by Craig, Guyenne and Kalisch \cite{CGK}, using
a Hamiltonian perturbation theory. It describes, as the Benjamin-Ono
equation, the evolution of weakly nonlinear dispersive internal long
waves at the interface of a two-layer system, one being infinitely
deep.

In \cite{CGK}, Craig, Guyenne and Kalisch (always using a
Hamiltonian perturbation theory) also derived a higher order
intermediate long wave equation.
\begin{equation}\label{ho.ILW}
\left\{\begin{array}[pos]{ll}
           \partial_tu-b\mathcal{F}_h
           \partial^2_xu+(a_1\mathcal{F}_h^2+a_2)\epsilon \partial_x^3u
           =cu\partial_xu-d\epsilon \partial_x(u\mathcal{F}_h
           \partial_xu+\mathcal{F}_h(u\partial_xu)) \\
           u(0)=\phi, \\
       \end{array} \right.
\end{equation}
where $\mathcal{F}_h$ is the Fourier multiplier $-i\coth(h\xi)$, $u$
is a real-valued solution, and $a_1, \ a_2, \ b, \ c, \ d$ and $h$
are positive constants. The same ill-posedness results also apply
for these equations.

These results are inspired by those from  Molinet, Saut and Tzvetkov
for the KPI equation \cite{MST1} and the Benjamin-Ono (and the ILW)
equation \cite{MST2}, (see also Bourgain \cite{Bou} and Tzvetkov
\cite{Tzv} for the KdV equation). It is worth notice that the
equation (\ref{ho.BO}) and the BO equation (as well as the equation
(\ref{ho.ILW}) and the ILW equation) share the same property of
ill-posedness of the flow in any Sobolev space $H^s(\R)$.

The rest of this paper is organized as follows: in Section 2, we
introduce a few notation, define the function spaces and state our
main results. In Section 3, we derive some linear estimates that we
use in Section 4 to prove our well-posedness results. Finally, in
Section 5, we deal with the ill-posedness results.

This work is part of my PhD Thesis at IMPA.\\

\noindent \textbf{Acknowledgments.} I would like to thank my advisor
Professor Felipe Linares for all his help and encouragements during
this work and Professor Jean-Claude Saut for helpful comments about
the ill-posedness part of this work. \\

\pagebreak
\section{Statements of the results}

\noindent \textbf{1. Some notations.} For any positive numbers $a$
and $b$, the notation $a \lesssim b$ means that there exists a
positive constant $c$ such that $a \le c b$. And we denote $a \sim
b$ when, $a \lesssim b$ and $b \lesssim a$.

Let $U_j(t)=e^{-t\partial^{2j+1}_x}$ be the unitary group (in
$H^s(\R)$) associated to the Airy equation, so that we have
\textit{via} Fourier transform
\begin{equation} \label{U_j}
U_j(t)\phi=\left(e^{(-1)^{j+1}i\xi^{2j+1}t}\widehat{\phi}\right)^{\vee},
\quad \forall \ t \in \R, \quad \forall \phi \in H^s(\R).
\end{equation}
The group $U_j$ commute with the operator of multiplication by $x$.
\begin{lemm} \label{lemm1}
Let $j \ge 1$ and $f \in \mathcal{S(\R)}$, then we have
\begin{equation} \label{lemm1.1}
xU_j(t)f=U_j(t)(xf)+(2j+1)tU_j(t)\partial_x^{2j}f \quad \forall \ t
\in \R.
\end{equation}
\end{lemm}
\noindent \textbf{Proof.} see \cite{Po}.\\
\linebreak \noindent
\textbf{2. Littlewood-Paley multipliers.} Throughout the paper, we
fix a cutoff function $\chi$ such that
\begin{equation} \label{chi}
\chi \in C_0^{\infty}(\R), \quad 0 \le \chi \le 1, \quad
\chi_{|_{[-1,1]}}=1 \quad \mbox{and} \quad  \mbox{supp}(\chi)
\subset [-2,2].
\end{equation}
We define
\begin{equation} \label{psi.1}
\psi(\xi):=\chi(\xi)-\chi(2\xi) \quad \mbox{and} \quad
\psi_l(\xi):=\psi(2^{-l}\xi),
\end{equation}
so that we have
\begin{equation} \label{psi.2}
\sum_{l \in \Z}\psi_l(\xi)=1, \ \forall \xi \neq 0 \quad \mbox{and}
\quad \mbox{supp}(\psi_l) \subset \{2^{l-1}\le |\xi| \le 2^{l+1}\}.
\end{equation}
Next, we define the Littlewood-Paley multipliers by
\begin{equation}  \label{Delta}
\Delta_lf=\left(\psi_l\widehat{f}\right)^{\vee}=(\psi_l)^{\vee}\ast
f \quad \forall f \in \mathcal{S}'(\R), \ \forall l \in \Z,
\end{equation}
and
\begin{equation}  \label{S}
S_lf=\sum_{k \le l}\Delta_kf  \quad \forall f \in \mathcal{S}'(\R),
\ \forall l \in \Z.
\end{equation}
More precisely we have that
\begin{equation}  \label{S}
S_0f=\left(\chi \widehat{f}\right)^{\vee}  \quad \forall f \in
\mathcal{S}'(\R),
\end{equation}
This means that $S_0$ is the operator of restriction in the low
frequencies. Note also that since
$(\psi_l)^{\vee}=2^l(\psi)^{\vee}(2^l\cdot)$,
$\|(\psi_l)^{\vee}\|_{L^1}=C$ and then, by Young's inequality we
have that for all $l \in \Z$
\begin{equation} \label{Delta2}
\|\Delta_lf\|_{L^p} \le C\|f\|_{L^p}, \ \forall \ f \in L^p, \
\forall p \in [1,+\infty].
\end{equation}
Combining this result with the integral Minkowski inequality, we
also deduce that
\begin{equation} \label{Delta2b}
\|\Delta_lf\|_{L^p_xL^q_t} \le C\|f\|_{L^p_xL^q_t}, \ \forall \ f
\in L^p_xL^q_t, \ \forall p, \ q \in [1,+\infty].
\end{equation}
We will need to commute $S_0$ and $\Delta_l$ with the operator of
multiplication by $x$
\begin{equation} \label{S_x}
[S_0,x]f= S_0^{'}f \quad \mbox{where} \quad S_0^{'}f=
\left((\frac{d}{d\xi}\chi)\widehat{f}\right)^{\vee}
\end{equation}
\begin{equation} \label{Delta_x}
[\Delta_l,x]f= \Delta_l^{'}f \quad \mbox{where} \quad
\Delta_l^{'}f=\left(2^{-l}(\frac{d}{d\xi}\psi)(2^{-l}\cdot)\widehat{f}\right)^{\vee}
\end{equation}
Finally, let $\tilde{\psi}$ be another smooth function supported in
$\R\setminus \{0\}$ such that $\tilde{\psi}=1$ on supp($\psi$). We
define $\tilde{\Delta}_l$ like $\Delta_l$ with $\tilde{\psi}$
instead of $\psi$ which yields in particular the following identity
\begin{equation} \label{Delta_tilde}
\tilde{\Delta}_l\Delta_l=\Delta_l.
\end{equation}

\noindent \textbf{3. Function spaces.} Let $1 \le p,q \le \infty$,
$T>0$,  the mixed \lq\lq space-time\rq\rq \ Lebesgue spaces are
defined by
$$L^p_xL^q_T := \{u:\R \times [-T,T] \rightarrow \R
\ \mbox{measurable} \ : \ \|u\|_{L^p_xL^q_T}<\infty \},$$ and
$$L^q_TL^p_x := \{u:\R \times [-T,T] \rightarrow \R
\ \mbox{measurable} \ : \ \|u\|_{L^q_TL^p_x}<\infty\},$$ where
\begin{equation} \label{LpLq}
\|u\|_{L^p_xL^q_T}:=\left(\int_{\R}\|u(x,\cdot)
\|_{L^q([-T,T])}^pdx\right)^{1/p},
\end{equation}
and
\begin{equation} \label{LqLp}
\|u\|_{L^q_TL^p_x}:=\left(\int_{-T}^T\|u(\cdot,t)
\|_{L^p(\R)}^qdt\right)^{1/q}.
\end{equation}
Next we derive the following Berstein's inequalities.
\begin{lemm} \label{lemm1b}
Let $f:\R\times [0,T]\rightarrow \C$ a smooth function and $p, \ q
\in [1, +\infty]$, then we have for all $j, \ l \in \N$,
\begin{equation} \label{lemm1b.1}
\|\Delta_l\partial_x^jf\|_{L^p_xL^q_T} \lesssim
2^{jl}\|\Delta_lf\|_{L^p_xL^q_T},
\end{equation}
and
\begin{equation} \label{lemm1b.2}
\|x\Delta_l\partial_x^jf\|_{L^p_xL^q_T} \lesssim
2^{jl}\|x\Delta_lf\|_{L^p_xL^q_T}+2^{j(l-1)}\|\Delta_lf\|_{L^p_xL^q_T}.
\end{equation}
\end{lemm}
\noindent \textbf{Proof.} We deduce from (\ref{Delta_tilde}),
(\ref{LpLq}), integral Minkowski's and Young's inequalities that
{\setlength\arraycolsep{2pt}
\begin{eqnarray*}
\|\Delta_l\partial_x^jf\|_{L^p_xL^q_T} &=&
\|\tilde{\Delta}_l\Delta_l\partial_x^jf\|_{L^p_xL^q_T}=
\|\|\partial_x^j(\tilde\psi_l)^{\vee}\ast
\Delta_lf(\cdot,t)\|_{L^q_T}\|_{L^p_x} \\
&\le&\|\int_{\R}|\partial_x^j(\tilde\psi_l)^{\vee}(y)| \|
\Delta_lf(x-y,t)\|_{L^q_T}dy\|_{L^p_x} \\
&\lesssim& \|\partial_x^j(\tilde\psi_l)^{\vee}\|_{L^1_x}
\|\Delta_lf\|_{L^p_xL^q_T}.
\end{eqnarray*}}
This implies (\ref{lemm1b.1}), since
$\|\partial_x^j(\tilde\psi_l)^{\vee}\|_{L^1_x}=c2^{lj}$. By a
similar argument, we have {\setlength\arraycolsep{2pt}
\begin{eqnarray*}
\|x\Delta_l\partial_x^jf\|_{L^p_xL^q_T} &\le&
\|x\left(|\partial_x^j(\tilde\psi_l)^{\vee}| \ast \|
\Delta_lf(\cdot,t)\|_{L^q_T}(x)\right)\|_{L^p_x} \\
&\lesssim& \|\partial_x^j(\tilde\psi_l)^{\vee}\|_{L^1_x}
\|x\Delta_lf\|_{L^p_xL^q_T}+\|x\partial_x^j(\tilde\psi_l)^{\vee}\|_{L^1_x}
\|\Delta_lf\|_{L^p_xL^q_T},
\end{eqnarray*}}
which implies inequality  (\ref{lemm1b.2}), since
$\|\partial_x^j(\tilde\psi_l)^{\vee}\|_{L^1_x}=c2^{lj}$ and
$\|x\partial_x^j(\tilde\psi_l)^{\vee}\|_{L^1_x}=c2^{(l-1)j}$.\hfill
$\square$

We will also use the fractional Sobolev spaces. Let $s \in \R$, then
$$ H^s(\R) := \{f \in \mathcal{S}'(\R) \ : \
(1+\xi^2)^{\frac{s}{2}}\widehat{f}(\xi) \in L^2(\R)\}$$ with the
norm
\begin{equation} \label {H}
\|f\|_{H^s}:=\|(1+\xi^2)^{s/2}\widehat{f}(\xi)\|_{L^2}.
\end{equation}
When $s=k \in \N$, it is well known (see for example \cite{Ste})
that
$$H^k(\R)=\{f \in L^2(\R) \ : \ \partial_x^lf \in L^2(\R), \ l= 0, 1
\cdots k\},$$ with the equivalent norm
\begin{equation} \label{L2.k}
\|f\|_{L^2_k}:=\sum_{l=0}^k\|\partial_x^lf\|_{L^2} \sim \|f\|_{H^k}.
\end{equation}
Similarly, it is possible to define weighted Sobolev spaces. Let $k
\in \N$, then
$$H^k(\R;x^2dx):=\{f \in L^2(\R;x^2dx) \ : \ \partial_x^lf \in
L^2(\R;x^2dx), \ l= 0, 1 \cdots k\},$$ with the norm
\begin{equation} \label{Hk.x}
\|f\|_{H^k(x^2dx)}:=\sum_{l=0}^k\|x\partial_x^lf\|_{L^2}.
\end{equation}

Finally, we recall the definition of the Besov spaces and define
weighted Besov spaces. Let $s \in \R$, $p,q \ge 1$, the non
homogeneous Besov space $\mathcal{B}_{p}^{s,q}(\R)$ is the
completion of the Schwartz space $\mathcal{S}(\R)$ under the norm
\begin{equation} \label{B}
\|f\|_{\mathcal{B}_{p}^{s,q}}:=\|S_0f\|_{L^p}+\|\{2^{ls}
\|\Delta_lf\|_{L^p}\}_{l \ge 0}\|_{l^q(\N)}.
\end{equation}
This definition naturally extends (even if $s \in \R$) for weighted
spaces. Let $s \in \R$, $p, \ q \ge 1$, then
$\mathcal{B}_{p}^{s,q}(\R;x^pdx)$ is the completion of the Schwartz
space $\mathcal{S}(\R)$ under the norm
\begin{equation} \label{B.x}
\|f\|_{\mathcal{B}_{p}^{s,q}(x^pdx)}:=\|xS_0f\|_{L^p}+
\|\{2^{ls}\|x\Delta_lf\|_{L^p}\}_{l \ge 0}\|_{l^q(\N)}.
\end{equation}
It is well known (see \cite{Tri}) that for all $s \in \R$
\begin{equation} \label{H-B}
H^s(\R)=\mathcal{B}_{2}^{s,2}(\R) \quad \mbox{and that} \quad
\|f\|_{H^s} \sim \|f\|_{\mathcal{B}_{2}^{s,2}}.
\end{equation}
Next we derive a similar result for weighted spaces in the case $s=k
\in \N$.
\begin{lemm} \label{lemm2}
Let $k \in \N$, $k \ge 1$ and $f \in \mathcal{S}(\R)$, then
\begin{equation} \label{lemm2.1}
\|f\|_{H^k(x^2dx)}+\|f\|_{H^{k-1}} \sim
\|f\|_{\mathcal{B}_{2}^{k,2}(x^2dx)}+\|f\|_{H^{k-1}}.
\end{equation}
\end{lemm}
\noindent \textbf{Proof.} We apply (\ref{S_x}), (\ref{Delta_x}),
(\ref{H-B}), the Plancherel theorem and the fact that the supports
of $\frac{d}{d\xi}\psi(2^{-l}\xi)$ are almost disjoint to get
{\setlength\arraycolsep{2pt}
\begin{eqnarray*}
\lefteqn{\|f\|_{\mathcal{B}_{2}^{k,2}(x^2dx)} =
\|xS_0f\|_{L^2}+\left(\sum_{l \ge 0}4^{lk}
\|x\Delta_lf\|_{L^2}^2\right)^{1/2}} \nonumber \\
&&\le   \|S_0(xf)\|_{L^2}+\|S_0'f\|_{L^2}+ \left(\sum_{l \ge
0}4^{lk}(\|\Delta_l(xf)\|_{L^2}+
\|\Delta_l'f\|_{L^2})^2\right)^{1/2} \nonumber \\
&&\lesssim \|xf\|_{\mathcal{B}_{2}^{k,2}}+
\left(\int_{\R}|(\frac{d}{d\xi}\chi)(\xi)\widehat{f}(\xi)|^2d\xi
+\sum_{l\ge0}\int_{\R}4^{l(k-1)}|(\frac{d}{d\xi}\psi)(2^{-l}\xi)
\widehat{f}(\xi)|^2d\xi\right)^{1/2} \\
&& \lesssim \|xf\|_{H^k}+\|\partial_x^{k-1}f\|_{L^2}.
\end{eqnarray*}}
Then we use (\ref{L2.k}) and the identity
$$\partial_x^l(xf)=l\partial_x^{l-1}f+x\partial_x^lf, \ \forall l \ge 1$$ to obtain that
\begin{equation} \label{lemm0.2.2}
\|f\|_{\mathcal{B}_{2}^{k,2}(x^2dx)} \lesssim
\|f\|_{H^k(x^2dx)}+\|f\|_{H^{k-1}}.
\end{equation}
The other inequality of (\ref{lemm2.1}) follows exactly by the same
argument. \hfill $\square$ \linebreak

\noindent \textbf{4. Statements of the results.}

\begin{theo} \label{theo1}
There exists $\delta>0$ such that for all $u_0 \in
\mathcal{B}^{2j+1/4,1}_2(\R) \cap \mathcal{B}^{1/4,1}_2(\R;x^2dx)$
with
\begin{equation} \label{theo1.1}
\beta=\|u_0\|_{\mathcal{B}_{2}^{2j+9/4,1}}+\|u_0\|_{\mathcal{B}_{2}^{1/4,1}(x^2dx)}
 \le \delta,
\end{equation}
there exists $T=T(\beta)$ such that $T(\beta)\nearrow +\infty$ when
$\beta \rightarrow 0$, a space $X_T$ such that
\begin{equation} \label{theo1.2}
X_T \hookrightarrow C([-T,T];\mathcal{B}^{2j+1/4,1}_2(\R) \cap
\mathcal{B}^{1/4,1}_2(\R;x^2dx))
\end{equation}
and a unique solution $u$ of (\ref{sKdV}) in $X_T$. Moreover, the
flow map is smooth from $\mathcal{B}^{2j+1/4,1}_2(\R) \cap
\mathcal{B}^{1/4,1}_2(\R;x^2dx)$ to $X_T$ near the origin.
\end{theo}

\begin{theo} \label{theo2}
Let $s>2j+1/4$, then there exists $\delta>0$ such that for all $u_0
\in H^s(\R) \cap \mathcal{B}^{s-2j,2}_2(\R;x^2dx)$ with
\begin{equation} \label{theo2.1}
\beta=\|u_0\|_{\mathcal{B}_{2}^{2j+1/4,1}}+\|u_0\|_{\mathcal{B}_{2}^{1/4,1}
(x^2dx)} \le \delta,
\end{equation}
there exists $T=T(\beta)$ such that $T(\beta)\nearrow +\infty$ when
$\beta \rightarrow 0$, a space $Y_{T,s}$ such that
\begin{equation} \label{theo2.2}
Y_{T,s} \hookrightarrow C([-T,T];H^s(\R) \cap
\mathcal{B}^{s-2j,2}_2(\R;x^2dx))
\end{equation}
and a unique solution $u$ of (\ref{sKdV}) in $Y_{T,s}$. Moreover,
the flow map is smooth from $H^s(\R) \cap
\mathcal{B}^{s-2j,2}_2(\R;x^2dx)$ to $Y_{T,s}$ near the origin.
\end{theo}

\begin{coro} \label{coro1}
Let $k \in \N$ such that $k>2j+1/4$, then the IVP (\ref{sKdV}) is
locally well-posed in the space $H^k(\R) \cap H^{k-2j}(\R;x^2dx)$
for small initial data.
\end{coro}

\noindent \textbf{Proof.} We know by Lemma \ref{lemm2}, that
$H^k(\R) \cap \mathcal{B}^{k-2j,2}_2(\R;x^2dx)=H^k(\R) \cap
H^{k-2j}(\R;x^2dx)$, then Corollary \ref{coro1} follows directly
from Theorem \ref{theo2}. \hfill $\square$

\begin{rema} \label{rem1}
Corollary \ref{coro1} improves the previous results in \cite{Arg}
for the non dissipative KdV-KV equation.
\end{rema}
 \pagebreak
Moreover, we have the following ill-posedness results for the IVP
(\ref{sKdV}).

\begin{theo} \label{theo3}
Let $s \in \R$ and $T>0$, suppose that there exists $k>j$ such that
$a_{0,k}\neq 0$, then, there does not exist any space $X_T$ such
that $X_T$ is continuously embedded in $C([-T,T];H^s(\R))$,
\textit{i.e.}
\begin{equation} \label{theo3.1}
\|u\|_{C([-T,T];H^s)} \lesssim \|u\|_{X_T}, \quad \forall \ u \in
X_T,
\end{equation}
and such that
\begin{equation} \label{theo3.2}
\|U_j(t)\phi\|_{X_T} \lesssim \|\phi\|_{H^s}, \quad \forall \ \phi
\in H^s(\R),
\end{equation}
and, for all $u$, $v \in X_ T$,
\begin{equation} \label{theo3.3}
\|\int_ 0^tU_j(t-t')\sum_{0\le l_1 \le l_2 \le 2j}
          a_{l_1,l_2}\partial_x^{l_1}u(t') \partial_x^{l_2}v(t')dt'\|_{X_T}
\lesssim \|u\|_{X_ T}\|v\|_{X_ T}.
\end{equation}
\end{theo}

\begin{theo} \label{theo4}
Let $s \in \R$, suppose that there exists $k>j$ such that
$a_{0,k}\neq 0$. Then, if the Cauchy problem (\ref{sKdV}) is locally
well-posed in $H^s(\R)$, the flow map data-solution
\begin{equation}\label{theo4.1}
S(t) :H^{s}(\R) \longrightarrow H^s(\R), \quad \phi \longmapsto u(t)
\end{equation}
is not $C^2$ at zero.
\end{theo}

These ill-posedness results also apply for the higher-order
equations (\ref{ho.BO}) and (\ref{ho.ILW})

\begin{theo} \label{theo5}
Let $s \in \R$. If the Cauchy problems (\ref{ho.BO}) and
respectively (\ref{ho.ILW}) are locally well-posed in $H^s(\R)$,
then the associated flow maps data-solution
\begin{equation}\label{theo5.1}
S^{hoBO}(t) :H^{s}(\R) \longrightarrow H^s(\R), \quad \phi
\longmapsto u(t),
\end{equation}
and respectively
\begin{equation}\label{theo5.2}
S^{hoILW}(t) :H^{s}(\R) \longrightarrow H^s(\R), \quad \phi
\longmapsto u(t)
\end{equation}
are not $C^2$ at zero.
\end{theo}

\pagebreak

\section{Linear estimates}

\noindent \textbf{1. Linear estimates for the free and the non
homogeneous evolutions.}
\begin{prop}[Kato type smoothing effect.]
Let $j \ge 1$. If $u_0 \in L^2(\R)$, then
\begin{equation} \label{Kse.1}
\|\partial^j_xU_j(t)u_0\|_{L^{\infty}_xL^2_t} \lesssim
\|u_0\|_{L^2}.
\end{equation}
Let $T>0$, then if $f \in L^1_xL^2_T$
\begin{equation} \label{Kse.2}
\|\int_0^t\partial^{j}_xU_j(t-t')f(\cdot,t')dt'\|_{L^{\infty}_TL^2_x}
\lesssim \|f\|_{L^1_xL^2_T},
\end{equation}
and
\begin{equation} \label{Kse.3}
\|\int_0^t\partial^{2j}_xU_j(t-t')f(\cdot,t')dt'\|_{L^{\infty}_xL^2_T}
\lesssim \|f\|_{L^1_xL^2_T}.
\end{equation}
\end{prop}
\noindent \textbf{Proof.} See \cite{KPV4}.

\begin{prop}[Maximal function estimate.] \label{MF}
If $u_0 \in \mathcal{S(\R)}$, then
\begin{equation} \label{mfe.1}
\|U_j(t)u_0\|_{L^{4}_xL^{\infty}_t} \lesssim \|D^{1/4}_xu_0\|_{L^2},
\end{equation}
and
\begin{equation} \label{mfe.2}
\|U_j(t)u_0\|_{L^{1}_xL^{\infty}_T} \lesssim
\|D^{1/4}_xu_0\|_{L^2}+\|D^{1/4}_x(xu_0)\|_{L^2}+T\|D^{1/4}_x\partial_x^{2j}u_0\|_{L^2}.
\end{equation}
\end{prop}
\noindent \textbf{Proof.} The estimate (\ref{mfe.1}) is due to
Kenig, Ponce and Vega \cite{KPV1} (see also the work of Kenig and
Ruiz \cite{KR} in the case $j=1$). We will prove the estimate
(\ref{mfe.2}) using (\ref{lemm1.1}), (\ref{mfe.1}) and H\"older's
inequality {\setlength\arraycolsep{2pt}
\begin{eqnarray*}
\|U_j(t)u_0\|_{L^{1}_xL^{\infty}_T} &=& \int_{|x| \le 1}
\sup_{[-T,T]}|U_j(t)u_0(x)|dx+\int_{|x| > 1}\frac{1}{|x|}\sup_{[-T,T]}|xU_j(t)u_0(x)|dx  \\
&\lesssim&
\|U_j(t)u_0\|_{L^{4}_xL^{\infty}_T}+\|U_j(t)(xu_0)\|_{L^{4}_xL^{\infty}_T}
+T\|U_j(t)\partial_x^{2j}u_0\|_{L^{4}_xL^{\infty}_T}  \\
&\lesssim&
\|D^{1/4}_xu_0\|_{L^2}+\|D^{1/4}_x(xu_0)\|_{L^2}+T\|D^{1/4}_x\partial_x^{2j}u_0\|_{L^2}.
\end{eqnarray*}}
\hfill $\square$ \linebreak

\begin{rema}
It is interesting to observe that the restriction on the $s$ in
Theorem \ref{theo2} ($s>2j+1/4$) appears in the estimate
(\ref{mfe.2}).
\end{rema}

\noindent \textbf{2. Linear estimates for phase localized
functions.} Following the ideas in \cite{MR2}, we  will derive
linear estimates for the phase localized free and nonhomogeneous
evolutions.

\begin{prop} \label{LE1}
Let $u_0 \in \mathcal{S}(\R)$, then we have for all $l \ge 0$
\begin{equation} \label{LE1.1}
\|\Delta_lU_j(t)u_0\|_{L^{\infty}_TL^{2}_x}=\|\Delta_lu_0\|_{L^{2}_x},
\end{equation}
and
\begin{equation} \label{LE1.2}
\|x\Delta_lU_j(t)u_0\|_{L^{\infty}_TL^{2}_x}\lesssim
\|x\Delta_lu_0\|_{L^{2}_x}+T2^{2jl}\|\Delta_lu_0\|_{L^{2}_x}.
\end{equation}
If $f:\R\times [0,T]\rightarrow \C$ is smooth, then we have for all
$l \ge 0$
\begin{equation} \label{LE1.3}
\|\int_0^t\Delta_lU_j(t-t')f(\cdot,t')dt'\|_{L^{\infty}_TL^{2}_x}\lesssim
2^{-jl}\|\Delta_lf\|_{L^1_xL^2_T},
\end{equation}
and
\begin{equation} \label{LE1.4}
\|\int_0^tx\Delta_lU_j(t-t')f(\cdot,t')dt'\|_{L^{\infty}_TL^{2}_x}\lesssim
2^{-jl}\|x\Delta_lf\|_{L^1_xL^2_T}
+T2^{jl}\|\Delta_lf\|_{L^1_xL^2_T}.
\end{equation}
\end{prop}

\noindent \textbf{Proof.} The identity (\ref{LE1.1}) follows
directly from the fact that $U_j$ is a unitary group in $L^2(\R)$.
To prove the estimate (\ref{LE1.2}), we will use (\ref{lemm1.1}),
(\ref{LE1.1}) and Plancherel's theorem
\begin{eqnarray*}
\|x\Delta_lU_j(t)u_0\|_{L^{\infty}_TL^{2}_x} &\le&
\|U_j(t)(x\Delta_lu_0)\|_{L^{\infty}_TL^{2}_x}
+(2j+1)T\|U_j(t)\partial^{2j}_x\Delta_lu_0\|_{L^{\infty}_TL^{2}_x} \\
&\lesssim&
\|x\Delta_lu_0\|_{L^{\infty}_TL^{2}_x}+T2^{2jl}\|\Delta_lu_0\|
_{L^{\infty}_TL^{2}_x}.
\end{eqnarray*}
The estimate (\ref{LE1.3}) follows from (\ref{Kse.2}), Plancherel's
theorem and the fact that $\Delta_l$ localize the frequency near
$|\xi|\sim2^l$. Next, we will prove the estimate (\ref{LE1.4}). The
identity (\ref{lemm1.1}), the estimate (\ref{Kse.2}) and the fact
the operator $x\Delta_l$ still localizes the frequency near $|\xi|
\sim 2^l$ (see the commutator identity (\ref{Delta_x}), imply that
\begin{eqnarray}
\lefteqn{\|\int_0^tx\Delta_lU_j(t-t')f(\cdot,t')dt'\|
_{L^{\infty}_TL^{2}_x}} \nonumber \\
 &&  \quad \quad \quad \quad \quad  \lesssim
\|\int_0^tU_j(t-t')(x\Delta_lf(\cdot,t'))dt'\|_{L^{\infty}_TL^{2}_x}
\nonumber \\ && \quad \quad \quad \quad \quad \quad
+T\|\int_0^t\Delta_lU_j(t-t')\partial_x^{2j}f(\cdot,t')dt'\|_{L^{\infty}_TL^{2}_x}
\nonumber \\ &&  \quad \quad \quad \quad \quad  \lesssim
2^{-jl}\|x\Delta_lf\|_{L^1_xL^2_T}+T2^{jl}\|\Delta_lf\|_{L^1_xL^2_T}.
\label{propLE1.1}
\end{eqnarray}
\hfill $\square$

\begin{prop} \label{LE2}
Let $u_0 \in \mathcal{S}(\R)$, then we have for all $l \ge 0$
\begin{equation} \label{LE2.1}
\|\Delta_lU_j(t)u_0\|_{L^{\infty}_xL^{2}_T}\lesssim
2^{-jl}\|\Delta_lu_0\|_{L^{2}_x},
\end{equation}
and
\begin{equation} \label{LE2.2}
\|x\Delta_lU_j(t)u_0\|_{L^{\infty}_xL^{2}_T}\lesssim
2^{-jl}\|x\Delta_lu_0\|_{L^{2}_x}+T2^{jl}\|\Delta_lu_0\|_{L^{2}_x}.
\end{equation}
If $f:\R\times [0,T]\rightarrow \C$ is smooth, then we have for all
$l \ge 0$
\begin{equation} \label{LE2.3}
\|\int_0^t\Delta_lU_j(t-t')f(\cdot,t')dt'\|_{L^{\infty}_xL^{2}_T}\lesssim
2^{-2jl}\|\Delta_lf\|_{L^1_xL^2_T},
\end{equation}
and
\begin{equation} \label{LE2.4}
\|\int_0^tx\Delta_lU_j(t-t')f(\cdot,t')dt'\|_{L^{\infty}_xL^{2}_T}\lesssim
2^{-2jl}\|x\Delta_lf\|_{L^1_xL^2_T} +T\|\Delta_lf\|_{L^1_xL^2_T}.
\end{equation}
\end{prop}

\noindent \textbf{Proof.} The proof is the same as for the
Proposition \ref{LE1} where we use (\ref{Kse.1}) and (\ref{Kse.3})
instead of (\ref{Kse.2}). \hfill $\square$

In order to derive a non homogeneous estimate for the localized
maximal function, we need the following lemma due to Molinet and
Ribaud (see \cite{MR2}) and inspired by a previous result of Christ
and Kiselev (see \cite{CK}).

\begin{lemm} \label{CK}
Let $L$ be a linear operator defined on space-time functions
$f(x,t)$ by
\begin{displaymath}
Lf(t)=\int_0^TK(t,t')f(t')dt',
\end{displaymath}
where $K:\mathcal{S}(\R^2)\rightarrow C(\R^3)$ and such that
\begin{displaymath}
\|Lf\|_{L^{p_1}_xL^{\infty}_T} \le C\|f\|_{L^{p_2}_xL^{q_2}_T},
\end{displaymath}
with $p_2,q_2<\infty$. Then,
\begin{displaymath}
\|\int_0^tK(t,t')f(t')dt'\|_{L^p_{x}L^{\infty}_T} \le
C\|f\|_{L^{p_2}_xL^{q_2}_T}.
\end{displaymath}
\end{lemm}

\begin{prop} \label{LE3}
Let $u_0 \in \mathcal{S}(\R)$, then we have for all $l \ge 0$
\begin{equation} \label{LE3.1}
\|\Delta_lU_j(t)u_0\|_{L^{1}_xL^{\infty}_T} \lesssim
2^{(\frac{1}{4}+2j)l}(1+T)\|\Delta_lu_0\|_{L^{2}_x}+2^{\frac{1}{4}l}
\|x\Delta_lu_0\|_{L^{2}_x}.
\end{equation}
If $f:\R\times [0,T]\rightarrow \C$ is smooth, then we have for all
$l \ge 0$
\begin{eqnarray} \label{LE3.2}
\lefteqn{\|\int_0^t\Delta_lU_j(t-t')f(\cdot,t')dt'\|_{L^{1}_xL^{\infty}_T}}
\nonumber \\ && \quad \quad \quad \quad \quad \quad \lesssim
2^{(\frac{1}{4}-j)l}\|x\Delta_lf\|_{L^1_xL^2_T}+(1+T)2^{(\frac{1}{4}+j)l}
\|\Delta_lf\|_{L^1_xL^2_T}. \label{LE3.2}
\end{eqnarray}
\end{prop}

\noindent \textbf{Proof.} To obtain the estimate (\ref{LE3.1}), we
apply (\ref{mfe.2}) with $\Delta_lu_0$ instead of $u_0$, then we use
Plancherel's theorem and the fact that the operators $\Delta_l$ and
$x\Delta_l$ localize the frequency near $|\xi| \sim 2^l$.

In order to prove the estimate (\ref{LE3.2}), we first need to
derive a "nonretarded" $L^4$-maximal function estimate. Note first
that duality and (\ref{mfe.1}) imply that
\begin{equation} \label{LE3.3}
\|\int_0^T\Delta_lU_j(-t)f(\cdot,t)dt\|_{L^{2}_x} \lesssim
2^{\frac{1}{4}l}\|\Delta_lf\|_{L^{4/3}_xL^1_T}.
\end{equation}
Then, we deduce combining (\ref{LE1.3}), (\ref{LE3.3}) and the
Cauchy-Schwarz inequality that for all $g \in L^{4/3}_xL^1_T$
\begin{eqnarray*}
\lefteqn{\int_{\R \times [0,T]}\left(\int_0^T\Delta_l
U_j(t-t')f(\cdot,t')dt'\right)g(x,t)dxdt} \\
&& =  \int_{\R}\left(\int_0^TU_j(-t')\Delta_lf(\cdot,t')dt'\right)
\left(\int_0^TU_j(-t)\tilde{\Delta}_l\overline{g}(\cdot,t)dt\right)dx \\
&& \le \|\int_0^TU_j(-t')\Delta_lf(\cdot,t')dt'\|_{L^2_x}
 \|\int_0^TU_j(-t)\tilde{\Delta}_l\overline{g}(\cdot,t)dt\|_{L^2_x}  \\
&& \lesssim
2^{-jl}\|\Delta_lf\|_{L^1_xL^{2}_T}2^{\frac{1}{4}l}\|g\|_{L^{4/3}_xL^{1}_T},
\end{eqnarray*}
so that by duality
\begin{equation} \label{LE3.4}
\|\int_0^T\Delta_lU_j(t-t')f(\cdot,t')dt'\|_{L^4_xL^{\infty}_T}
\lesssim 2^{(\frac{1}{4}-j)l}\|\Delta_lf\|_{L^1_xL^{2}_T}.
\end{equation}
Then, we use Lemma \ref{CK} to obtain the corresponding "retarded"
estimate
\begin{equation} \label{LE3.5}
\|\int_0^t\Delta_lU_j(t-t')f(\cdot,t')dt'\|_{L^4_xL^{\infty}_T}
\lesssim 2^{(\frac{1}{4}-j)l}\|\Delta_lf\|_{L^1_xL^{2}_T}.
\end{equation}

We are now able to derive the $L^1_xL^{\infty}_T$ estimate for the
non homogeneous term. We have by H\"older's inequality
\begin{eqnarray}
\lefteqn{\|\int_0^t\Delta_lU_j(t-t')f(\cdot,t')dt'\|_{L^1_xL^{\infty}_T}}
\nonumber \\  &&  = \int_{|x|\le1}\sup_{t
\in[-T,T]}\left|\int_0^t\Delta_l
U_j(t-t')f(\cdot,t')dt'\right|dx \nonumber  \\
&& \quad  +\int_{|x|>1}\frac{1}{|x|}\sup_{t\in[-T,T]}
\left|\int_0^tx\Delta_lU_j(t-t')f(\cdot,t')dt'\right|dx   \label{LE3.6} \\
&& \lesssim
\|\int_0^t\Delta_lU_j(t-t')f(\cdot,t')dt'\|_{L^4_xL^{\infty}_T}
+\|\int_0^tx\Delta_lU_j(t-t')f(\cdot,t')dt'\|_{L^4_xL^{\infty}_T}
\nonumber
\end{eqnarray}
Thus, we deduce from (\ref{lemm1.1}), (\ref{LE3.5}) and
(\ref{LE3.6}) that
\begin{eqnarray}
\lefteqn{\|\int_0^t\Delta_lU_j(t-t')f(\cdot,t')dt'\|_{L^1_xL^{\infty}_T}}
\nonumber\\
&&  \lesssim
(2^{(\frac{1}{4}-j)l}+T2^{(\frac{1}{4}+j)l})\|\Delta_l\|_{L^1_xL^{2}_T}
+2^{(\frac{1}{4}-j)l}\|x\Delta_lf\|_{L^1_xL^{2}_T}, \label{LE3.7}
\end{eqnarray}
which leads to (\ref{LE3.2}), since $l \ge 0$. \hfill $\square$

\begin{rema}
All the results in Propositions \ref{LE1}, \ref{LE2} and \ref{LE3}
are still valid with $S_0$ instead of $\Delta_l$ and $l=0$.
\end{rema}

\section{Proof of Theorems \ref{theo1} and \ref{theo2}}

\noindent \textbf{Proof of Theorem \ref{theo1}.} Consider the
integral equation associated to (\ref{sKdV})
\begin{equation} \label{IE}
u(t)=F(u)(t),
\end{equation}
where
\begin{equation} \label{IEb}
F(u)(t):=U_j(t)u_0+\int_0^tU_j(t-t')\sum_{0\le j_1 + j_2 \le 2j}
a_{j_1,j_2}\partial_x^{j_1}u(t')
\partial_x^{j_2}u(t')dt'.
\end{equation}
Let $T>0$, define the following seminorms:
{\setlength\arraycolsep{2pt}
\begin{eqnarray}
N_1^T(u)
&=&\|S_0u\|_{L^{\infty}_TL^2_{x}}+\sum_{l=1}^{\infty}2^{(2j+\frac{1}{4})l}
\|\Delta_lu\|_{L^{\infty}_TL^2_{x}}, \label{theo1.3} \\
N_2^T(u)&=&\|xS_0u\|_{L^{\infty}_TL^2_{x}}+\sum_{l=1}^{\infty}
2^{\frac{1}{4}l}\|x\Delta_lu\|_{L^{\infty}_TL^2_{x}}, \label{theo1.4} \\
P_1^T(u)&=&\|S_0u\|_{L^{\infty}_xL^2_{T}}+\sum_{l=1}^{\infty}
2^{(3j+\frac{1}{4})l}\|\Delta_lu\|_{L^{\infty}_xL^2_{T}}, \label{theo1.5} \\
P_2^T(u)&=&\|xS_0u\|_{L^{\infty}_xL^2_{T}}+\sum_{l=1}^{\infty}
2^{(j+\frac{1}{4})l}\|x\Delta_lu\|_{L^{\infty}_xL^2_{T}}, \label{theo1.6} \\
M^T(u)&=&\|S_0u\|_{L^1_xL^{\infty}_{T}}+\sum_{l=1}^{\infty}
\|\Delta_lu\|_{L^1_xL^{\infty}_{T}}. \label{theo1.7}
\end{eqnarray}}
Then, we define the Banach space
\begin{equation} \label{theo1.8}
X_T=\{u \in C([-T,T];\mathcal{B}^{2j+1/4,1}_2(\R) \cap
\mathcal{B}^{1/4,1}_2(\R;x^2dx)) \ : \ \|u\|_{X_T}<\infty\},
\end{equation}
where
\begin{equation} \label{theo1.9}
\|u\|_{X_T}=N_1^T(u)+N_2^T(u)+P_1^T(u)+P_2^T(u)+M^T(u).
\end{equation}

We deduce from (\ref{LE1.1}), (\ref{LE1.2}), (\ref{LE2.1}),
(\ref{LE2.2}) and (\ref{LE3.1}) that
\begin{equation} \label{theo1.10}
\|U_j(t)u_0\|_{X_T} \lesssim (1+T)
\left(\|u_0\|_{\mathcal{B}^{2j+1/4,1}_2}+\|u_0\|_{\mathcal{B}^{1/4,1}_2(x^2dx)}\right),
\end{equation}
and from (\ref{LE1.3}), (\ref{LE1.4}), (\ref{LE2.3}), (\ref{LE2.4})
and (\ref{LE3.2}) that
\begin{eqnarray} \label{theo1.11}
\lefteqn{\|\int_0^tU_j(t-t')\sum_{0\le j_1 + j_2 \le 2j}
a_{j_1,j_2}\partial_x^{j_1}u(t')\partial_x^{j_2}v(t')dt'\|_{X_T}}
\nonumber \\ && \lesssim (1+T)\sum_{0\le j_1 + j_2 \le 2j}
|a_{j_1,j_2}|
\left(\|S_0(\partial_x^{j_1}u\partial_x^{j_2}v)\|_{L^1_xL^2_T}
\right. \nonumber \\ && \quad
+\sum_{l=1}^{\infty}2^{(j+\frac{1}{4})l}
\|\Delta_l(\partial_x^{j_1}u\partial_x^{j_2}v)\|_{L^1_xL^2_T}
+\|xS_0(\partial_x^{j_1}u\partial_x^{j_2}v)\|_{L^1_xL^2_T}\nonumber
\\ && \quad
\left.+\sum_{l=1}^{\infty}2^{(\frac{1}{4}-j)l}
\|x\Delta_l(\partial_x^{j_1}u\partial_x^{j_2}v)\|_{L^1_xL^2_T}\right).
\end{eqnarray}

In order to estimate the nonlinear term
$\sum_{l=1}^{\infty}2^{(j+\frac{1}{4})l}\|\Delta_l(\partial^{j_1}_xu
\partial^{j_2}_xv)\|_{L^1_xL^2_T}$,
we observe that {\setlength\arraycolsep{2pt}
\begin{eqnarray}
\Delta_l(fg) &=& \Delta_l\left((S_0f+\sum_{r \ge
1}\Delta_rf)(S_0g+\sum_{k \ge 1}\Delta_kg)\right)
\nonumber \\
&=& \Delta_l\left(S_0fS_0g+\sum_{r \ge 1}\Delta_rfS_rg+\sum_{r \ge
1}\Delta_rgS_{r-1}f\right), \label{theo1.12}
\end{eqnarray}}
where $f=\partial^{j_1}_xu$ and $g=\partial^{j_2}_xv$. First, since
$\Delta_l(S_0uS_0v)=0$ for $l\ge 3$ and since the operators
$\Delta_l$ are uniformly bounded (in $l$) in $L^1$, we have by
H\"older's inequality
\begin{equation} \label{theo1.13}
\sum_{l\ge 1}2^{(j+\frac{1}{4})l}\|\Delta_l(S_0\partial^{j_1}_xuS_0
\partial^{j_2}_xv)\|_{L^1_xL^2_T}
\lesssim \|S_0u\|_{L^{\infty}_xL^2_T}\|S_0v\|_{L^1_xL^{\infty}_T}.
\end{equation}
In order to estimate the second term on the right-hand side of
(\ref{theo1.12}), we notice, since the term
$\Delta_r\partial^{j_1}_xuS_{r}\partial^{j_2}_xv$ is localized in
frequency in the set $|\xi|\le2^{r+3}$ and the operator $\Delta_l$
only sees the frequency in the set $2^{l-1}\le|\xi|\le2^{l+1}$, that
\begin{equation} \label{theo1.14}
\Delta_l\left(\sum_{r=1}^{\infty}\Delta_rfS_{r}g \right)=
\Delta_l\left(\sum_{r\ge l-3}\Delta_rfS_{r}g \right).
\end{equation}
Then, we only have to estimate terms of the form $\Delta_l(\sum_{r
\ge l}\Delta_r\partial_x^{j_1}uS_{r}\partial_x^{j_2}v)$. Using
H\"older's inequality, the estimate (\ref{lemm1b.1}), and the fact
that {\setlength\arraycolsep{2pt}
\begin{eqnarray}
\|S_r\partial_x^{j_2}v\|_{L^1_xL^{\infty}_T} &\le&
\|S_0v\|_{L^1_xL^{\infty}_T}+\sum_{k=1}^r
\|\Delta_k\partial_x^{j_2}v\|_{L^1_xL^{\infty}_T}
\nonumber \\
&\le&\|S_0v\|_{L^1_xL^{\infty}_T}+\sum_{k=1}^r2^{j_2k}
\|\Delta_kv\|_{L^1_xL^{\infty}_T} \lesssim 2^{j_2r}M^T(v),
\label{theo1.15b}
\end{eqnarray}}
we deduce that
\begin{eqnarray}
\lefteqn{\sum_{l \ge 1}2^{(j+\frac{1}{4})l}\|\Delta_l(\sum_{r \ge
l}\Delta_r\partial_x^{j_1}uS_{r}\partial_x^{j_2}v)\|_{L^1_xL^2_T}}
\nonumber \\ && \quad \quad \quad \le \sum_{l \ge
1}2^{(j+\frac{1}{4})l}\sum_{r \ge
l}\|\Delta_r\partial_x^{j_1}u\|_{L^{\infty}_xL^2_T}
\|S_r\partial_x^{j_2}v\|_{L^1_xL^{\infty}_T} \nonumber \\ && \quad
\quad \quad \le
M^T(v)\sum_{r\ge1}\left(\sum_{l=1}^r2^{(j+\frac{1}{4})l}\right)
2^{(j_1+j_2)r}\|\Delta_ru\|_{L^{\infty}_xL^2_T}
 \nonumber \\ && \quad \quad \quad  \lesssim
M^T(v)P_1^T(u) \le \|u\|_{X_T}\|v\|_{X_T}, \label{theo1.15}
\end{eqnarray}
since $j_1+j_2\le2j$. Thus, we obtain, gathering (\ref{theo1.12}),
(\ref{theo1.13}), (\ref{theo1.14}) and (\ref{theo1.15}) that
\begin{equation} \label{theo1.16}
\sum_{l=1}^{\infty}2^{(j+\frac{1}{4})l}\|\Delta_l(\partial_x^{j_1}
u\partial_x^{j_2}v)\|_{L^1_xL^2_T}\lesssim \|u\|_{X_T}\|v\|_{X_T}.
\end{equation}
We apply exactly the same strategy to estimate the other nonlinear
term $\sum_{l=1}^{\infty}2^{(\frac{1}{4}-j)l}\|x\Delta_l
(\partial_x^{j_1}u\partial_x^{j_2}v)\|_{L^1_xL^2_T}$. Then, we have
only to estimate terms of the form
$$\sum_{l=1}^{\infty}2^{(\frac{1}{4}-j)l}\|x\Delta_l(\sum_{r\ge
l}\Delta_r\partial_x^{j_1}uS_r\partial_x^{j_2}v)\|_{L^1_xL^2_T}.$$
For this, we combine the same ideas as for the estimate
(\ref{theo1.16}) with the commutator identity (\ref{Delta_x}) and
the fact that the operators $\Delta_l'$ are also uniformly bounded
(in $l$) in $L^1$ to deduce that {\setlength\arraycolsep{2pt}
\begin{eqnarray} \sum_{l \ge
1}2^{(\frac{1}{4}-j)l}\|x\Delta_l(\partial_x^{j_1}u
\partial_x^{j_2}v)\|_{L^1_xL^2_T}
&\lesssim&  \sum_{l \ge 1}2^{(\frac{1}{4}-j)l}\sum_{r \ge
l}\|x\Delta_r\partial_x^{j_1}uS_r\partial_x^{j_2}v\|_{L^{\infty}_xL^2_T}
\nonumber \\ &+& \sum_{l \ge 1}2^{-(\frac{3}{4}+j)l}\sum_{r \ge
l}\|\Delta_r\partial_x^{j_1}uS_r\partial_x^{j_1}v\|_{L^{\infty}_xL^2_T} \nonumber \\
&\lesssim&
 M^T(v)(P_1^T(u)+P^T_2(u)). \label{theo1.17}
\end{eqnarray}}
Thus, we deduce from (\ref{theo1.11}), (\ref{theo1.16}) and
(\ref{theo1.17}) that
\begin{eqnarray} \label{theo1.18}
\lefteqn{\|\int_0^tU_j(t-t')\sum_{ j_1 + j_2 \le 2j}
a_{j_1,j_2}\partial_x^{j_1}u(t')
\partial_x^{j_2}u(t')dt'\|_{X_T}} \nonumber \\ &&
\quad \quad  \quad \quad \quad \quad \quad \quad \quad \quad \quad
\quad \quad \quad \lesssim (1+T)\|u\|_{X_T}\|v\|_{X_T}.
\end{eqnarray}
Then, we use (\ref{theo1.10}) and  (\ref{theo1.18}) to deduce that
there exists a constant $C>0$ such that
\begin{equation} \label{theo1.19}
\|F(u)\|_{X_T} \le C(1+T)\left(\|u_0\|_{\mathcal{B}^{9/4,1}_2}
+\|u_0\|_{\mathcal{B}^{1/4,1}_2(x^2dx)}+\|u\|_{X_T}^2\right), \quad
\forall \ u \in X_T,
\end{equation}
and
\begin{equation} \label{theo1.20}
\|F(u)-F(v)\|_{X_T} \le C(1+T)(\|u\|_{X_T}+ \|v\|_{X_T})
\|u-v\|_{X_T}, \quad \forall \ u, v \in X_T.
\end{equation}

Let $X_T(a):=\{u\in X_T \ : \ \|u\|_{X_T}\le a\}$ the closed ball of
$X_T$ with radius $a$. $X_T(a)$ equipped with the metric induced by
the norm $\|\cdot\|_{X_T}$ is a complete metric space. If we choose
\begin{equation} \label{theo1.21}
\beta=\|u_0\|_{\mathcal{B}^{9/4,1}_2}+\|u_0\|_{\mathcal{B}^{1/4,1}_2
(x^2dx)}\le \delta< \min\{(\frac{1}{4C})^2,1\},
\end{equation}
\begin{equation} \label{theo1.22}
a=\sqrt{\beta}, \quad \mbox{and} \quad T=\frac{1}{4C\sqrt{\beta}},
\end{equation}
we have that
\begin{equation} \label{theo1.22b}
2C(1+T)a<1.
\end{equation}
Then, we deduce from (\ref{theo1.19}) and (\ref{theo1.20}) that the
operator $F$ is a contraction in $X_T(a)$ (up to the persistence
property) and so, by the Picard fixed point theorem, there exists a
unique solution of (\ref{IE}) in $X_T(a)$.

The proof of persistence, uniqueness and smoothness of the map
follows by standard arguments (see for example \cite{KPV1b}). \hfill
$\square$

\noindent \textbf{Proof of Theorem 2.}

\begin{lemm} \label{lemma4}
Let $s>2j+1/4$, then the injection
\begin{equation} \label{lemma4.1}
 H^s(\R) \cap \mathcal{B}^{s-2j,2}_2(\R;x^2dx)\hookrightarrow
 \mathcal{B}^{2j+1/4,1}_2(\R)\cap \mathcal{B}^{1/4,1}_2(\R;x^2dx)
\end{equation}
is continuous.
\end{lemm}

\noindent \textbf{Proof.} Let $s>2j+1/4$ and $f \in H^s(\R)$. We
obtain using the Cauchy-Schwarz inequality that
{\setlength\arraycolsep{2pt}
\begin{eqnarray}
\|f\|_{\mathcal{B}^{2j+1/4,1}_2}&=&\|S_0f\|_{L^2}
+\sum_{l \ge 0}2^{ls}\|\Delta_lf\|_{L^2}2^{l(2j+1/4-s)} \nonumber \\
&\le& \|S_0f\|_{L^2}+\left(\sum_{l \ge
1}4^{(2j+1/4-s)l}\right)^{1/2}
\left(\sum_{l \ge 1}4^{sl}\|\Delta_lf\|^2_{L^2}\right)^{1/2} \nonumber \\
&\lesssim& \|f\|_{\mathcal{B}^{s,2}_2}\sim \|f\|_{H^s}.
\label{lemma4.2}
\end{eqnarray}}
Similarly, we get
\begin{equation} \label{lemma4.3}
\|f\|_{\mathcal{B}^{1/4,1}_2(x^2dx)} \lesssim
\|f\|_{\mathcal{B}^{s-2j,2}_2(x^2dx)},
\end{equation}
when $s > 2j+1/4$ and then, (\ref{lemma4.2}) and (\ref{lemma4.3})
yield (\ref{lemma4.1}). \hfill $\square$

Now, let $s>2j+1/4$. Exactly as in the proof of Theorem 1, we want
to apply a fixed point theorem to solve the integral equation
(\ref{IE}) in some good function space. In this way, define the
following semi-norm
\begin{equation} \label{theo2.3}
\|u\|_{X_{T,s}}=N_{1,s}^T(u)+N_{2,s}^T(u)+P_{1,s}^T(u)+P_{2,s}^T(u),
\end{equation}
where {\setlength\arraycolsep{2pt}
\begin{eqnarray}
N_{1,s}^T(u)
&=&\|S_0u\|_{L^{\infty}_TL^2_{x}}+\left(\sum_{l=1}^{\infty}
4^{sl}\|\Delta_lu\|^2_{L^{\infty}_TL^2_{x}}\right)^{1/2}, \label{theo2.4} \\
N_{2,s}^T(u)
&=&\|xS_0u\|_{L^{\infty}_TL^2_{x}}+\left(\sum_{l=1}^{\infty}
4^{(s-2j)l}\|x\Delta_ju\|^2_{L^{\infty}_TL^2_{x}}\right)^{1/2}, \label{theo2.5} \\
P_{1,s}^T(u)&=&\|S_0u\|_{L^{\infty}_xL^2_{T}}+\left(\sum_{l=1}^{\infty}
4^{(s+j)l}\|\Delta_lu\|^2_{L^{\infty}_xL^2_{T}}\right)^{1/2}, \label{theo2.6} \\
P_{2,s}^T(u)&=&\|xS_0u\|_{L^{\infty}_xL^2_{T}}+\left(\sum_{l=1}^{\infty}
4^{(s-j)l}\|x\Delta_lu\|^2_{L^{\infty}_xL^2_{T}}\right)^{1/2}.
\label{theo2.7}
\end{eqnarray}}
If $u_0 \in H^s(\R) \cap \mathcal{B}^{s-2j,2}_2(\R;x^2dx)$, by Lemma
\ref{lemma4}, it makes sense to define
\begin{equation} \label{theo2.8}
\lambda_s :=
\frac{\|u_0\|_{\mathcal{B}^{2j+1/4,1}_2}+\|u_0\|_{\mathcal{B}^{1/4,1}_2(x^2dx)}}
{\|u_0\|_{H^s}+\|u_0\|_{\mathcal{B}^{s-2j,2}_2(x^2dx)}}.
\end{equation}
Then, let $Y_{T,s}$ be the Banach space
\begin{equation} \label{theo2.9}
Y_{T,s}=\{u \in C([-T,T];H^s(\R) \cap
\mathcal{B}^{s-2j,2}_2(\R;x^2dx)) \ \mbox{such that} \
\|u\|_{Y_{T,s}}<\infty\},
\end{equation}
where
\begin{equation} \label{theo2.10}
 \|u\|_{Y_{T,s}}=\|u\|_{X_T}+\lambda_s\|u\|_{X_{T,s}}.
\end{equation}

We deduce from (\ref{LE1.1}), (\ref{LE1.2}), (\ref{LE2.1}),
(\ref{LE2.2}), (\ref{LE3.1}), (\ref{theo2.8}) and (\ref{theo2.10})
that
\begin{equation} \label{theo2.11}
\|U_j(t)u_0\|_{Y_{T,s}} \lesssim (1+T)
\left(\|u_0\|_{\mathcal{B}^{2j+1/4,1}_2}+\|u_0\|_{\mathcal{B}^{1/4,1}_2(x^2dx)}\right).
\end{equation}
In order to estimate the nonlinear term of (\ref{IEb}) in the norm
$\|\cdot\|_{Y_{T,s}}$, we remember (\ref{theo1.18}), and then it
only remains to derive estimates of the form
\begin{equation} \label{theo2.12}
\|\int_0^tU_j(t-t'))(\partial_x^{j_1}u(t')
\partial_x^{j_2}u(t'))dt'\|_{X_{T,s}} \lesssim
(1+T)\|u\|_{Y_{T,s}}\|v\|_{Y_{T,s}}.
\end{equation}
In this way, we use (\ref{LE1.3}), (\ref{LE1.4}), (\ref{LE2.3}),
(\ref{LE2.4}) and (\ref{LE3.2}) and argue as in the proof of Theorem
\ref{theo1} to estimate the left-hand side of (\ref{theo2.12}) by
some terms of the form
\begin{equation} \label{theo2.14}
A= \left(\sum_{l \ge
1}4^{(s-j)l}\|\Delta_l(\sum_{r=l}^{\infty}\Delta_r\partial_x^{j_1}u
S_r\partial_x^{j_2}v)\|^2_{L^1_xL^2_{T}}\right)^{1/2},
\end{equation}
\begin{equation} \label{theo2.15}
B= \left(\sum_{l \ge
1}4^{(s-3j)l}\|x\Delta_l(\sum_{r=l}^{\infty}\Delta_r\partial_x^{j_1}u
S_r\partial_x^{j_2}v)\|^2_{L^1_xL^2_{T}}\right)^{1/2},
\end{equation}
and some others harmless terms. We next estimate $A$, we get from
(\ref{Delta2b}), H\"older's inequality, (\ref{lemm1b.1}) and
(\ref{theo1.15b}), the inequality
\begin{equation} \label{theo2.16}
A \le M^T(v)\left(\sum_{l \ge
0}4^{(s+j)l}\left(\sum_{r=l}^{\infty}\|\Delta_ru
\|_{L^{\infty}_xL^2_{T}}\right)^2\right)^{1/2}.
\end{equation}
Then, define
\begin{equation} \label{theo2.17}
\gamma_r=2^{(s+j)l}\|\Delta_ru\|_{L^{\infty}_xL^2_{T}} \quad
\mbox{and note that} \quad \|\{\gamma_r\}_r\|_{l^2(\N)} \le
P^T_{1,s}(u).
\end{equation}
We deduce by (\ref{theo2.16}), a change of index and Minkowski's
inequality that {\setlength\arraycolsep{2pt}
\begin{eqnarray*}
A &\le&
M^T(v)\|\{\sum_{r=l}^{\infty}2^{(s+j)(l-r)}\gamma_r\}_l\|_{l^2(\N)}
= M^T(v)\|\{\sum_{k \ge 0}2^{-(s+j)k}\gamma_{l+k}\}_l\|_{l^2(\N)} \\
&\le& M^T(v)\sum_{k \ge
0}2^{-(s+j)k}\|\{\gamma_{l+k}\}_l\|_{l^2(\N)} \le
M^T(v)\|\{\gamma_{l}\}_l\|_{l^2(\N)}\sum_{k \ge 0}2^{-(s+j)k},
\end{eqnarray*}}
then (\ref{theo2.17}) implies that
\begin{equation} \label{theo2.18}
A \lesssim  P^T_{1,s}(u)M^T(v).
\end{equation}
Analogously, we obtain a similar estimate for $B$
\begin{equation} \label{theo2.19}
B \lesssim  P^T_{2,s}(u)M^T(v).
\end{equation}
Thus, (\ref{theo2.18}) and (\ref{theo2.19}) yield (\ref{theo2.12})
and we conclude the proof of Theorem \ref{theo2} as for Theorem
\ref{theo1} using (\ref{theo2.11}) and (\ref{theo2.12}) instead of
(\ref{theo1.10}) and (\ref{theo1.18}). \hfill $\square$

\begin{rema}
Unless we can use the strategy of Kenig, Ponce and Vega in
\cite{KPV3} to show well-posedness for the IVPs (\ref{ho.BO}) and
(\ref{ho.ILW}) in weighted Sobolev spaces, it is not clear wether
the technique used here applies or not.
\end{rema}

\section{Ill-posedness results}

In the proofs of Theorems \ref{theo3} and \ref{theo4}, we will
suppose, for simplicity, that the nonlinearity $\sum_{0\le l_1 \le
l_2 \le 2j}a_{l_1,l_2}\partial_x^{l_1}u \partial_x^{l_2}u$ has
the form $\partial_x^k(u^2)$ with $k > j$. \\

\noindent \textbf{Proof of Theorem \ref{theo3}.} The key point of
the proof is the following algebraic relation
\begin{lemm} \label{lemm6}
Let $j \in \N$ such that $j \ge 1$ and $\xi, \xi_1 \in \R$, then
\begin{equation} \label{lemm6.1}
\xi_1^{2j+1}+(\xi-\xi_1)^{2j+1}-\xi^{2j+1}=(\xi-\xi_1)Q_{2j}(\xi,\xi_1),
\end{equation}
where
\begin{equation} \label{lemm6.2}
Q_{2j}(\xi,\xi_1)=\sum_{l=0}^{2j}((-1)^lC^l_{2j}-1)\xi^{2j-l}\xi_1^l
\end{equation}
and $C^l_{n}=\frac{n!}{l!(n-l)!}$.
\end{lemm}
Note that $\xi-\xi_1$ does not divide $Q_{2j}(\xi,\xi_1)$.

Let $s \in \R$, $k, j \in \N$ such that $k>j$ and $T>0$. Suppose
that there exists a space $X_T$ such as in Theorem \ref{theo3}. Take
$\phi$, $\psi \in H^s(\R)$, and define $u(t)=U_j(t)\phi$ and
$v(t)=U_j(t)\psi$. Then, we use (\ref{theo3.1}), (\ref{theo3.2}) and
(\ref{theo3.3}) to deduce that
\begin{equation} \label{theo3.4}
\|\int_ 0^tU_j(t-t')\partial_x^k[(U_j(t')\phi)
(U_j(t')\psi)]dt'\|_{H^s} \lesssim \|\phi\|_{H^s}\|\psi\|_{H^s}.
\end{equation}
We will show that (\ref{theo3.4}) fails for an appropriate pair of
$\phi$, $\psi$, which  would lead to a contradiction.

Define $\phi$ and $\psi$ by
\begin{equation} \label{theo3.5}
\phi=(\alpha^{-1/2}\chi_{I_ 1})^{\vee}
\end{equation}
and
\begin{equation} \label{theo3.5b}
\psi=(\alpha^{-1/2}N^{-s}\chi_{I_ 2})^{\vee}
\end{equation}
where
\begin{equation} \label{theo3.6}
N \gg 1, \quad  0<\alpha\ll1, \quad I_1=[\alpha/2,\alpha] \quad
\mbox{and} \quad  I_2=[N,N+\alpha]
\end{equation}
Note first that
\begin{equation} \label{theo3.6}
\|\phi\|_{H^s}\sim\|\psi\|_{H^s}\sim 1.
\end{equation}
Then, we use the algebraic relation (\ref{lemm6.1}), the definition
of the unitary group $U_j$ and the definition of $\phi$ and $\psi$
to estimate the Fourier transform of the left-hand side  of
(\ref{theo1.4})
\begin{eqnarray} \label{theo3.7}
\lefteqn{\left(\int_ 0^tU_j(t-t')\partial_x^k[(U_j(t')\phi )
(U_j(t')\psi)]dt'\right)^{\wedge}(\xi)} \nonumber\\
& & =\int_0^te^{(-1)^{j+1}i(t-t')\xi^{2j+1}}(i\xi)^k
(e^{(-1)^{j+1}it(\cdot)^{2j+1}}\widehat{\phi})\ast
(e^{(-1)^{j+1}it(\cdot)^{2j+1}}\widehat{\psi})(\xi)dt' \nonumber\\
& & =\int_{\R}e^{(-1)^{j+1}it\xi^{2j+1}}(i\xi)^k
\widehat{\psi}(\xi_1)\widehat{\phi}(\xi-\xi_1)
\int_0^te^{(-1)^{j+1}it'Q_{2j}(\xi,\xi_1)(\xi-\xi_1)}dt'd\xi_1 \nonumber\\
& & =
\int_{\R}e^{(-1)^{j+1}it\xi^{2j+1}}(i\xi)^k\widehat{\psi}(\xi_1)
\widehat{\phi}(\xi-\xi_1)
\frac{e^{(-1)^{j+1}it(\xi-\xi_1)Q_{2j}(\xi,\xi_1)}-1}{(-1)^{j+1}
i(\xi-\xi_1)Q_{2j}(\xi,\xi_1)}d\xi_1.
  \nonumber\\
& & \sim \frac{e^{(-1)^{j+1}it\xi^{2j+1}}\xi^k}{\alpha N^s} \int_{
\scriptsize{ \left\{  \begin{array}[pos]{ll}
                  \xi_1 \in I_2   \\          \xi-\xi_1 \in I_1  \\
                \end{array}           \right. }}
                \frac{e^{(-1)^{j+1}it(\xi-\xi_1)Q_{2j}(\xi,\xi_1)}-1}
                {(\xi-\xi_1)Q_{2j}(\xi,\xi_1)}d\xi_1.
\end{eqnarray}
When $\xi-\xi_1 \in I_1$ and $\xi_1 \in I_2$, we have that
$|(\xi-\xi_1)Q_{2j}(\xi,\xi_1)| \sim \alpha N^{2j}$. We choose
$\alpha = N^{-2j-\epsilon}$, with $0< \epsilon <1$ so that
\begin{equation} \label{theo3.8}
|(\xi-\xi_1)Q_{2j}(\xi,\xi_1)| \sim  N^{-\epsilon} \ll 1
\end{equation}
and
\begin{equation} \label{theo3.9}
\frac{e^{(-1)^{j+1}it(\xi-\xi_1)Q_{2j}(\xi,\xi_1)}-1}{(\xi-\xi_1)Q_{2j}
(\xi,\xi_1)}=ct+o(N^{-\epsilon})
\end{equation}
where $c \in \C$. We are now able to give a lower bound for the
left-hand side of (\ref{theo3.4})
\begin{equation} \label{theo3.10}
\|\int_ 0^tU_j(t-t')\partial_x^k[(U_j(t')\phi)
(U_j(t')\psi)]dt'\|_{H^s} \gtrsim
\frac{N^s}{N^s\alpha}N^k\alpha^{1/2}\alpha.
\end{equation}
Thus we conclude from (\ref{theo3.4}), (\ref{theo3.6}) and
(\ref{theo3.10}) that
\begin{equation} \label{theo3.11}
N^k\alpha^{1/2}=N^{k-j-\epsilon/2} \lesssim 1, \quad \forall \ N \gg
1,
\end{equation}
which is a contradiction since $k>j$. \hfill $\square$ \\

\begin{rema} Since the class of equation (\ref{sKdV}) often appears in
physical situations where the function $u$ is needed to be
real-valued, it is interesting to notice that Theorems \ref{theo3}
and \ref{theo4} are also valid if we ask the functions to be real.
Actually  take $\phi_1=\Re{\phi}$ and $\psi_1=\Re{\psi}$ instead of
$\phi$ and $\psi$, then
\begin{equation} \label{rema}
\widehat{\phi_1}=\frac{\alpha^{-1/2}}{2}\chi_{\{\alpha/2\le |\xi|
\le \alpha\}} \quad \mbox{and} \quad
\widehat{\psi_1}=\frac{\alpha^{-1/2}N^{-s}}{2} \chi_{\{N \le |\xi|
\le N+\alpha\}},
\end{equation}
and so we can conclude the proof as above.
\end{rema}

\noindent \textbf{Proof of Theorem \ref{theo4}.} Let $s \in \R$ and
$k, j \in \N$ such that $k>j$. Suppose that there exists $T>0$ such
that the Cauchy problem (\ref{sKdV}) is locally well-posed in
$H^s(\R)$ in the time interval $[0,T]$ and that its flow map
solution $S^{j,k} :H^s(\R)\longrightarrow C([0,T];H^s(\R))$ is $C^2$
at the origin. When $\phi \in H^s(\R)$, we will denote
$u_{\phi}(t)=S^{j,k}(t)\phi$ the solution of the Cauchy problem
(\ref{sKdV}) with initial data $\phi$. This means that $u_{\phi}$ is
a solution of the integral equation
\begin{equation} \label{theo4.2b}
u(t):=U_j(t)u_0+\int_0^tU_j(t-t')\partial_x^k(u^2)dt'.
\end{equation}
When $\phi$ and $\psi$ are in $H^s(\R)$, we use the fact that the
nonlinearity $\partial_x^k(uv)$ is a bilinear symmetric application
to compute the Fr\'echet derivative of $S^{j,k}(t)$ at $\psi$ in the
direction $\phi$
\begin{equation} \label{theo4.2}
d_{\psi}S^{j,k}(t)\phi=U_j(t)\phi+2\int_0^tU_j(t-t')\partial_x^k
(u_{\psi}(t')d_{\psi}S^{j,k}(t')\phi)dt'.
\end{equation}
Since the Cauchy problem (\ref{sKdV}) is supposed to be
well-posed,we know  using the uniqueness that
$S^{j,k}(t)0=u_{0}(t)=0$ and then we deduce from ({\ref{theo4.2}})
that
\begin{equation} \label{theo4.3}
d_{0}S^{j,k}(t)\phi=U_j(t)\phi.
\end{equation}
Using (\ref{theo4.2}), we compute the second Fr\'echet derivative at
the origin in the direction $(\phi,\psi)$
{\setlength\arraycolsep{2pt}
\begin{eqnarray*}
d_0^2S^{j,k}(t)(\phi,\psi) &=& d_0(d \
S^{j,k}(t)\phi)\psi=\frac{\partial}{\partial_{\beta}}
(\beta \mapsto d_{\beta \psi}S^{j,k}(t)\phi)_{|_{\beta=0}} \\
&=&  2\int_0^tU_j(t-t')\partial_x^k(d_{\beta \psi}S^{j,k}(t')\psi
d_{\beta \psi}S^{j,k}(t')\phi)dt'_{|_{\beta=0}} \\
&&  +2\int_0^tU_j(t-t')\partial_x^k(u_{\beta \psi}(t')d_{\beta
\psi}^2 S^{j,k}(t')(\phi,\psi))dt'_{|_{\beta=0}}.
\end{eqnarray*}}
Thus we deduce using (\ref{theo4.3}) that
\begin{equation} \label{theo4.4}
d_0^2S^{j,k}(t)(\phi,\psi)=2\int_0^tU_j(t-t')\partial_x^k[(U_j(t')\psi)(U_j(t')\phi)]dt'.
\end{equation}
The assumption of $C^2$ regularity of $S^{j,k}(t)$ at the origin
would imply that $d_0^2S^{j,k}(t) \in \mathcal{B}(H^s(\R)\times
H^s(\R),H^s(\R))$, which would lead to the following inequality
\begin{equation} \label{theo4.5}
\|d_0^2S^{j,k}(t)(\phi,\psi)\|_{H^s(\R)} \lesssim \|\phi\|_{H^s(\R)}
\|\psi\|_{H^s(\R)}, \quad \forall \ \phi, \ \psi \in H^s(\R).
\end{equation}
But (\ref{theo4.5}) is equivalent to (\ref{theo3.4}) which has been
shown
to fail in the proof of Theorem \ref{theo3}. \hfill $\square$ \\

\noindent \textbf{The case of the higher-order Benjamin-Ono and
intermediate long wave equations.} In order to study the Cauchy
problems (\ref{ho.BO}) (respectively (\ref{ho.ILW})), we define
$V_1$ (respectively $V_2(t)$) the unitary group in $H^s(\R)$
associated to the linear part of the equations, \textit{i.e.}
\begin{equation} \label{V}
V_k(t)\phi=\left(e^{ip_k(\xi)t}\widehat{\phi}\right)^{\vee}, \quad
k=1,2, \quad \forall \ t \in \R, \quad \forall \phi \in H^s(\R),
\end{equation}
where
\begin{displaymath} \label{p_1}
p_1(\xi)=b|\xi|\xi+a\epsilon\xi^3,
\end{displaymath}
and
\begin{displaymath} \label{p_1}
p_2(\xi)=b\coth(h\xi)\xi^2+(a_1\coth^2(h\xi)+a_2)\epsilon\xi^3.
\end{displaymath}
We denote by $f_1$ (respectively $f_2$) the nonlinearity of the
equations (\ref{ho.BO}) (respectively (\ref{ho.ILW})), \textit{i.e.}
\begin{displaymath}
f_1(u)=cu\partial_xu-d\epsilon
\partial_x(uH\partial_xu+H(u\partial_xu)),
\end{displaymath}
and
\begin{displaymath}
f_2(u)=cu\partial_xu-d\epsilon
\partial_x(u\mathcal{F}_h\partial_xu+\mathcal{F}_h(u\partial_xu)).
\end{displaymath}
Then, we have the analogous of Theorem \ref{theo3} for the equations
(\ref{ho.BO}) and (\ref{ho.ILW}).

\begin{theo} \label{theo6}
Let $s \in \R$, $T>0$ and $k \in \{1,2\}$. Then, there does not
exist any space $X_T$ such that $X_T$ is continuously embedded in
$C([-T,T];H^s(\R))$, \textit{i.e.},
\begin{equation} \label{theo6.1}
\|u\|_{C([-T,T];H^s)} \lesssim \|u\|_{X_T}, \quad \forall \ u \in
X_T,
\end{equation}
and such that
\begin{equation} \label{theo6.2}
\|V_k(t)\phi\|_{X_T} \lesssim \|\phi\|_{H^s}, \quad \forall \ \phi
\in H^s(\R),
\end{equation}
and
\begin{equation} \label{theo6.3}
\|\int_ 0^tV_k(t-t')f_k(u)(t')dt'\|_{X_T} \lesssim \|u\|^2_{X_ T},
\quad \forall \ u \in X_ T.
\end{equation}
\end{theo}
Theorem \ref{theo5} is a consequence of Theorem \ref{theo6} (see the
proof
of Theorem \ref{theo4}). \\

\noindent \textbf{Proof of Theorem \ref{theo6}.} Let $s \in \R$,
$T>0$ and $k \in \{1,2\}$. Suppose that there exists a space $X_T$
such as in Theorem \ref{theo6}. Take $\phi \in H^s(\R)$, and define
$u(t)=V_k(t)\phi$. Then, we use (\ref{theo6.1}), (\ref{theo6.2}) and
(\ref{theo6.3}) to see that
\begin{equation} \label{theo6.4}
\|\int_ 0^tV_k(t-t')f_k(V_k(t'))\phi)dt'\|_{H^s} \lesssim
\|\phi\|^2_{H^s}.
\end{equation}
We will show that (\ref{theo6.4}) fails for an appropriate choice of
$\phi$, which would lead to a contradiction.

Define $\phi$ by \footnote{We can also take $\Re{\phi}$ instead of
$\phi$ (see the remark after the proof of Theorem \ref{theo3}).}
\begin{equation} \label{theo6.6}
\phi=\left(\alpha^{-1/2}\chi_{I_ 1}+\alpha^{-1/2}N^{-s}\chi_{I_
2}\right)^{\vee}
\end{equation}
where
\begin{equation} \label{theo6.7}
N \gg 1, \quad  0<\alpha\ll1, \quad I_1=[\alpha/2,\alpha] \quad
\mbox{and} \quad  I_2=[N,N+\alpha]
\end{equation}
Note first that
\begin{equation} \label{theo6.7b}
\|\phi\|_{H^s}\sim 1.
\end{equation}
Then, the same computation as for (\ref{theo6.7}) leads to
\begin{equation} \label{theo6.8}
\left(\int_ 0^tV_k(t-t')f_k((V_k(t')\phi)dt'\right)^{\wedge}(\xi)
\sim g_1(\xi,t)+g_2(\xi,t)+g_3(\xi,t),
\end{equation}
where,
\begin{displaymath}
g_1(\xi,t)= \frac{e^{itp(\xi)}}{\alpha} \int_{
\tiny{\begin{array}[pos]{ll}
                  \xi_1 \in I_1   \\          \xi-\xi_1 \in I_1  \\
                \end{array}}}
\tilde{f}_k(\xi,\xi_1)\frac{e^{it(p(\xi_1)+p(\xi-\xi_1)-p(\xi))}-1}
{i(p(\xi_1)+p(\xi-\xi_1)-p(\xi))}d\xi_1,
\end{displaymath}

\begin{displaymath}
g_2(\xi,t)= \frac{e^{itp(\xi)}}{\alpha N^{2s}} \int_{
\tiny{\begin{array}[pos]{ll}
                  \xi_1 \in I_2   \\          \xi-\xi_1 \in I_2  \\
                \end{array}}}
\tilde{f}_k(\xi,\xi_1)\frac{e^{it(p(\xi_1)+p(\xi-\xi_1)-p(\xi))}-1}
{i(p(\xi_1)+p(\xi-\xi_1)-p(\xi))}d\xi_1,
\end{displaymath}

{\setlength\arraycolsep{2pt}
\begin{eqnarray*}
g_3(\xi,t)&=& \frac{e^{itp(\xi)}}{\alpha N^{s}} \left( \int_{ \tiny{
\begin{array}[pos]{ll}
                  \xi_1 \in I_1   \\          \xi-\xi_1 \in I_2  \\
                \end{array}  }}
\tilde{f}_k(\xi,\xi_1)\frac{e^{it(p(\xi_1)+p(\xi-\xi_1)-p(\xi))}-1}
{i(p(\xi_1)+p(\xi-\xi_1)-p(\xi))}d\xi_1, \right.\\
&& + \left. \int_{ \tiny{\begin{array}[pos]{ll}
                  \xi_1 \in I_2   \\          \xi-\xi_1 \in I_1  \\
                \end{array}}}
\tilde{f}_k(\xi,\xi_1)\frac{e^{it(p(\xi_1)+p(\xi-\xi_1)-p(\xi))}-1}
{i(p(\xi_1)+p(\xi-\xi_1)-p(\xi))}d\xi_1\right),
\end{eqnarray*}}
and
\begin{displaymath}
\tilde{f}_1(\xi,\xi_1)=c\xi_1-d\epsilon(\xi|\xi_1|+|\xi|\xi_1),
\end{displaymath}
or
\begin{displaymath}
\tilde{f}_2(\xi,\xi_1)=c\xi_1-d\epsilon(\xi\coth(\xi_1)\xi_1+\coth(\xi)\xi\xi_1).
\end{displaymath}
Since the supports of $g_1(\cdot,t)$, $g_2(\cdot,t)$ and
$g_3(\cdot,t)$ are disjoint, we use (\ref{theo6.8}) to bound by
below the left-hand side of (\ref{theo6.4})
\begin{equation} \label{theo6.9}
\|\int_ 0^tV_k(t-t')f_k((V_k(t'))\phi)dt'\|_{H^s} \ge \|(g_
3)^{\vee}(\xi,t)\|_ {H^s}.
\end{equation}

We notice that the function $p_k$ is smooth and that
\begin{equation} \label{theo6.10}
|p'_k(\xi)| \lesssim 1+|\xi|^2.
\end{equation}
Thus, when $\xi_1 \in I_1$ and $\xi-\xi_1 \in I_2$ or $\xi-\xi_1 \in
I_1$ and $\xi_1 \in I_2$, we have that $|\xi| \sim N$, and we use
(\ref{theo6.10}) and the mean value theorem to get the estimate
\begin{equation} \label{theo6.11}
|p(\xi_1)+p(\xi-\xi_1)-p(\xi)| \lesssim \alpha N^2.
\end{equation}
Hence we choose $\alpha = N^{-2-\epsilon}$, with $0< \epsilon <1$,
to get
\begin{equation} \label{theo6.12}
\left|\frac{e^{it(p(\xi_1)+p(\xi-\xi_1)-p(\xi))}-1}{p(\xi_1)+p(\xi-\xi_1)-p(\xi)}\right|
=|t|+o(N^{-\epsilon}).
\end{equation}
We are now able to give a lower bound for $\|(g_ 3)^{\vee}(\xi,t)\|_
{H^s}$
\begin{equation} \label{theo6.13}
\|(g_ 3)^{\vee}(\xi,t)\|_ {H^s} \gtrsim \frac{N^s}{N^s\alpha}
\left(N^2\alpha^{1/2}\alpha-N\alpha\alpha^{1/2}\alpha\right) \gtrsim
N^2\alpha^{1/2}.
\end{equation}
Thus, we conclude from (\ref{theo6.4}), (\ref{theo6.7b}),
(\ref{theo6.9}) and (\ref{theo6.13}) that
\begin{equation} \label{theo6.14}
N^2\alpha^{1/2}=N^{1-\epsilon/2} \lesssim 1, \quad \forall \ N \gg
1,
\end{equation}
which is a contradiction. \hfill $\square$

\noindent {\small{UFRJ, Institute of Mathematics, \\
P.O. Box 68530 - Cidade Universit\'aria. \\
Ilha do Fund\~ao. CEP 21945-970 \\
Rio de Janeiro, RJ, Brazil.}

\noindent E-mail: pilod@impa.br, didier@im.ufrj.br
\end{document}